\pdfoutput=1
\RequirePackage{ifpdf}
\ifpdf 
\documentclass[pdftex]{sigma}
\else
\documentclass{sigma}
\fi

\numberwithin{equation}{section}
\newtheorem{Theorem}{Theorem}[section]
\newtheorem{Corollary}[Theorem]{Corollary}
\newtheorem{Lemma}[Theorem]{Lemma}
\newtheorem{Proposition}[Theorem]{Proposition}
{\theoremstyle{definition}

\newtheorem{Example}[Theorem]{Example}
}

\begin{document}

\newcommand{\arXivNumber}{1311.3880}

\allowdisplaybreaks

\renewcommand{\thefootnote}{$\star$}

\renewcommand{\PaperNumber}{068}

\FirstPageHeading

\ShortArticleName{Groupoid Actions on Fractafolds}

\ArticleName{Groupoid Actions on Fractafolds\footnote{This paper is a~contribution to the Special Issue on
Noncommutative Geometry and Quantum Groups in honor of Marc A.~Rief\/fel.
The full collection is available at
\href{http://www.emis.de/journals/SIGMA/Rieffel.html}{http://www.emis.de/journals/SIGMA/Rieffel.html}}}

\Author{Marius IONESCU~$^\dag$ and Alex KUMJIAN~$^\ddag$}

\AuthorNameForHeading{M.~Ionescu and A.~Kumjian}

\Address{$^\dag$~Department of Mathematics, United States Naval Academy, Annapolis, MD, 21402-5002, USA}
\EmailD{\href{mailto:felijohn@gmail.com}{felijohn@gmail.com}}

\Address{$^\ddag$~Department of Mathematics, University of Nevada, Reno, NV, 89557, USA}
\EmailD{\href{mailto:alex@unr.edu}{alex@unr.edu}}

\ArticleDates{Received February 04, 2014, in f\/inal form June 21, 2014; Published online June 28, 2014}

\Abstract{We def\/ine a~bundle over a~totally disconnected set such that each f\/iber is homeo\-morphic to a~fractal blowup.
We prove that there is a~natural action of a~Renault--Deaconu groupoid on our fractafold bundle and that the resulting
action groupoid is a~Renault--Deaconu groupoid itself.
We also show that when the bundle is locally compact the associated $C^*$-algebra is primitive and has a~densely def\/ined
lower-semicontinuous trace.}


\Classification{28A80; 22A22; 46L55; 46L05}

\rightline{\emph{Dedicated to Marc A.~Rieffel on the occasion of his 75th birthday}}

\renewcommand{\thefootnote}{\arabic{footnote}}
\setcounter{footnote}{0}

\section{Introduction}

The goal of this paper is to f\/ind and analyze symmetries of fractals associated to iterated function systems
$(F_1,\dots, F_{N})$ and to study the associated $C^*$-algebras that might arise from the dynamics.
In~\cite{Str_CJM98} Stricharz constructed a~family of fractafold blowups of the invariant set of an iterated function
system which is parameterized by inf\/inite words in the alphabet $\{1, \dots, N \}$ and observed that two such blowups
are naturally homeomorphic if the parametrizing words are eventually the same.
We endow these fractafold blowups with the inductive limit topology and assemble them into a~fractafold bundle~$L$.

In general there do not appear to be any natural nontrivial symmetries of a~generic blowup but Stricharz's observation
suggests that we look for symmetries of the bundle instead.
Indeed we show that the homeomorphisms between f\/ibers observed by Stricharz give rise to a~natural groupoid action
on~$L$, the fractafold bundle.
This groupoid action and the associated action groupoid $\tilde{\mathcal{G}}$ constitute the main focus for this work.

We prove that the there is a~local homeomorphism $\tilde{\sigma}$ on~$L$ such that $\tilde{\mathcal{G}}$ is isomorphic
to the Renault--Deaconu groupoid associated to $\tilde{\sigma}$ and, in particular, $\tilde{\mathcal{G}}$ is \'etale.
We also prove that~$\tilde{\mathcal{G}}$ is topologically principal and has a~dense orbit.
If the iterated function system satisf\/ies the open set condition then we construct a~$\tilde{\mathcal{G}}$-invariant
measure on the unit space $\tilde{\mathcal{G}}^0$.

Now suppose that~$L$ is locally compact.
Then the associated $C^*$-algebra, $C^*(\tilde{\mathcal{G}})$, is primitive.
If, in addition the iterated function system satisf\/ies the open set condition, then the associated $C^*$-algebra has
a~densely def\/ined lower semi-continuous trace.

We begin by reviewing some of the background material and by proving some general results about Renault--Deaconu
groupoids in Section~\ref{sec:Groupoidactions}.
These results will be useful in our analysis of the fractafold bundle~$L$ and the associated action groupoid
$\tilde{\mathcal{G}}$ in Section~\ref{sec:groupoidaction}.
In the f\/inal section we consider some examples to illustrate the theory.
We show that the action groupoid is not in general minimal.
We also point out cases when the fractafold bundle~$L$ fails to be locally compact.

\section{Renault--Deaconu groupoid and groupoid actions}
\label{sec:Groupoidactions}

In this section we prove that if $\mathcal{G}$ is the Renault--Deaconu grou\-po\-id associated to a~local
homeomorphism~$\sigma$ on a~topological space~$X$ such that $\mathcal{G}$ acts on a~topological space~$Z$, then~$\sigma$
lifts to a~natural local homeomorphism $\tilde{\sigma}$ on~$Z$.
Moreover, we show that the resulting action groupoid is isomorphic to the Renault--Deaconu groupoid associated to
$\tilde{\sigma}$.
We begin by reviewing some of the background material that we need.
While some of the motivation for this work came from the theory of $C^*$-algebras associated to groupoids, many of our
results hold for topological groupoids that are not necessarily locally compact.
Since some of the groupoids in our examples are not locally compact, the only generic assumption on the topological
spaces and topological groupoids that we make is that they are Hausdorf\/f.

Let~$G$ be a~Hausdorf\/f topological groupoid (see~\cite{Ren_cartan_08}, cf.\ also~\cite{Ren_LNM793}).
Then the structure maps are continuous and, in addition, both the range map~$r$ (where $r(x)=xx^{-1}$) and the source
map~$s$ (where $s(x)=x^{-1}x$) are open.
We write $G^{0}$ for the \emph{unit space} of~$G$.
The groupoid~$G$ is said to be \emph{\'{e}tale} if~$s$ is a~local homeomorphism (or equivalently,~$r$ is a~local
homeomorphism).
A~subset $S \subset G$ is called a~\emph{$G$-set} or a~\emph{bisection} if the restrictions, $r|_S$, $s|_S$, are both
injective.
If~$G$ is \'{e}tale, then it has a~cover of open~$G$-sets and the restriction of either~$r$ or~$s$ to an open~$G$-set is
a~homeomorphism onto an open subset of $G^0$.

Let~$G$ be topological groupoid.
If the set of points in $G^0$ with trivial isotropy,
\begin{gather*}
\{x \in G^0 \mid s(\gamma) = r(\gamma) = x
\
\text{only if}
\
\gamma = x \},
\end{gather*}
is dense in $G^0$, we say that~$G$ is \emph{topologically principal} (see~\cite[Def\/inition~3.5(ii)]{Ren_cartan_08}).

If~$G$ is a~\emph{locally compact} \'{e}tale groupoid we let $C^*(G)$ denote the full $C^*$-algebra of~$G$ and
$C^*_r(G)$ denote the reduced $C^*$-algebra of~$G$.
If, in addition,~$G$ is amenable\footnote{See~\cite[Proposition 2.2.13]{anaren:amenable00} for dif\/ferent
characterizations of amenability for topological groupoids.} the canonical quotient map $C^*(G) \to C^*_r(G)$ is an
isomorphism.
We regard $C_0(G^0)$ as an abelian $C^*$-subalgebra of both.

Let~$X$ be a~topological space and let $\sigma:X\to X$ be a~local homeomorphism on~$X$.
The \emph{Renault--Deaconu groupoid} associated to~$\sigma$
\cite{ADC_BSMF97, Dea_TAMS95,De_Ku_Mu_JOT01,Ren_LNM793,Ren_CLA00} is
\begin{gather*}
\mathcal{G}= \mathcal{G}(X, \sigma) = \big\{(x,m-n,y)\in X\times \mathbb{Z}\times X:\sigma^m(x)=\sigma^n(y)\big\}.
\end{gather*}
Two elements $(x,m,y)$ and $(z,n,w)$ in $\mathcal{G}$ are composable if and only if $y=z$ and, in this case, their
product is
\begin{gather*}
(x,m,y)(y,n,w)=(x,m+n,w).
\end{gather*}
The inverse of an element in $\mathcal{G}$ is def\/ined~by
\begin{gather*}
(x,n,y)^{-1}=(y,-n,x).
\end{gather*}
Thus the range and source maps are given by $r(x,n,y)=(x,0,x)$ and $s(x,n,y)=(y,0,y)$.
Hence $\mathcal{G}^{0}$, the unit space of $\mathcal{G}$, may be identif\/ied with~$X$ via the map $(x,0,x)\mapsto x$ (in
the sequel we will often make this identif\/ication without comment).
A~basis for the topology consists of sets of the form
\begin{gather*}
\mathcal{G}(U,m,n,V)=\big\{(x,m-n,y):\sigma^m(x)=\sigma^n(y), x\in U, y\in V\big\},
\end{gather*}
where~$m$ and $n\in\mathbb{N}$,~$U$ and~$V$ are open subsets of~$X$ such that both $\sigma^m|_U$ and $\sigma^n|_V$ are
injective and $\sigma^m(U)=\sigma^n(V)$.
Note that the range map~$r$ induces a~homeomorphism $\mathcal{G}(U,m,n,V) \cong U$.
Hence, with this topology $\mathcal{G}$ is an \'{e}tale groupoid.
If~$X$ is locally compact then $\mathcal{G}=\mathcal{G}(X,\sigma)$ is a~locally compact groupoid.

Let~$G$ be a~topological groupoid with unit space $G^0 = X$.
Let~$Z$ be a~topological space and let $\rho: Z\to X$ be a~continuous open map.
A~continuous map $\alpha: G *Z \to Z$ (where $G *Z= \{(\gamma,z) \mid s(\gamma)=\rho(z) \}$), write $\alpha(\gamma, z) =
\gamma\cdot z$, such that $\rho(\gamma\cdot z)=r(\gamma)$, $\gamma_2\cdot(\gamma_1\cdot z)=(\gamma_2\gamma_1)\cdot z$,
and $\rho(z)\cdot z=z$, for all $\gamma,\gamma_1,\gamma_2\in G$ and $z\in Z$ is said to be a~\emph{left action} of~$G$
on~$Z$ (see~\cite{MuReWi_JOT87}).
Then $G*Z$ may itself be endowed with the structure of a~topological groupoid where the topology is inherited from $G
\times Z$; it is called the \emph{(left) action groupoid} (associated to~$\alpha$).
If~$G$ and~$Z$ are locally compact, then so is $G *Z$.
Recall that $(\gamma_1,z_1)$ and $(\gamma_2,z_2)\in G*Z$ are composable if $z_1=\gamma_2\cdot z_2$ and the product is
given by $(\gamma_1,z_1)(\gamma_2,z_2)=(\gamma_1\gamma_2,z_2)$.
The inverse of an element $(\gamma,z)$ is $(\gamma^{-1},\gamma\cdot z)$.
Therefore $(G *Z)^0 = \{(\rho(z), z) \mid z \in Z \}$ and thus the unit space may be identif\/ied with~$Z$ using the
projection onto the second factor $(\rho(z), z) \mapsto z$; note that $s(\gamma, z) = (s(\gamma), z)$ and $r(\gamma, z)
= (r(\gamma), \gamma\cdot z)$.
If~$V$ is an open subset of~$G$ and~$W$ is an open subset of~$Z$, then $V * W:= (V \times W) \cap (G *Z)$ is an open
subset of $G*Z$.

\begin{Lemma}
\label{lem:etale}
With notation as above let $\alpha: G *Z \to Z$ be a~left action and suppose that~$G$ is \'etale.
Then the action groupoid $G*Z$ is itself \'etale and the action map $\alpha: G *Z \to Z$ is open.
\end{Lemma}
\begin{proof}
Let $(\gamma_0, z_0) \in G *Z$ and let~$V$ be an open~$G$-set containing $\gamma_0$.
Then $s|_V$ yields a~ho\-meo\-morphism $V \cong s(V)$.
Set $W:= \rho^{-1}(s(V))$ and observe that $V * W$ is an open neighborhood of~$(\gamma_0, z_0)$.
Then since $s(\gamma, z) = (s(\gamma), z)$, we have $s(V * W) = s(V) * W = \rho(W) * W$ is open and the restriction
$s|_{V * W}$ is a~continuous injection; the inverse map $(\rho(z), z) \mapsto ((s|_V)^{-1}(\rho(z)), z)$ is continuous.
Hence, $s|_{V * W}$ is a~homeomorphism onto an open subset of $(G *Z)^0$.
So $s: G *Z \to (G *Z)^0$ is a~local homeomorphism and $G*Z$ is \'etale.

Let $\pi: G*Z \to Z$ be the projection map onto the second factor and observe that $\alpha = \pi\circ r$.
Then since $r(G*Z) = (G *Z)^0$ is open and both~$\pi$ and~$r$ are open maps, the action map $\alpha: G *Z \to Z$ is
open.
\end{proof}

We next show that if the groupoid acting is a~Renault--Deaconu groupoid the action groupoid is itself of the same type.

\begin{Theorem}
\label{thm:lift}
Let~$\sigma$ be a~local homeomorphism on a~topological space~$X$ and suppose that the Renault--Deaconu groupoid
$\mathcal{G}=\mathcal{G}(X,\sigma)$ acts on the left on the topological space~$Z$.
Define $\tilde{\sigma}:Z\to Z$ via
\begin{gather*}
\tilde{\sigma}(z)=(\sigma(\rho(z)),-1,\rho(z))\cdot z.
\end{gather*}
Then $\tilde{\sigma}$ is a~local homeomorphism on~$Z$ such that $\rho\circ\tilde{\sigma}=\sigma\circ\rho$.
Moreover, the left action groupoid $\mathcal{G}\ast Z$ is isomorphic to $\tilde{\mathcal{G}} = \mathcal{G}(Z,
\tilde{\sigma})$, the Renault--Deaconu groupoid associated to $\tilde{\sigma}$, via the map $\Psi:\tilde{\mathcal{G}}\to
\mathcal{G}\ast Z$ given~by
\begin{gather}
\label{eq:isomorphism}
\Psi((t,m-n,s))=((\rho(t),m-n,\rho(s)),s),
\end{gather}
where $(t,m-n,s) \in \tilde{\mathcal{G}} \subset Z \times \mathbb{Z}\times Z$.
\end{Theorem}
\begin{proof}
Let $z\in Z$.
We need to f\/ind an open neighborhood~$U$ of~$z$ such that $\tilde{\sigma}(U)$ is open and $\tilde{\sigma}|_U:U\to
\tilde{\sigma}(U)$ is a~homeomorphism.
Since~$\sigma$ is a~local homeomorphism there is an open neighborhood $V\subset X$ of $\rho(z)$ such that $\sigma(V)$ is
open and $\sigma|_V:V\to \sigma(V)$ is a~homeomorphism.
Let $U=\rho^{-1}(V)$.
We prove f\/irst that $\tilde{\sigma}|_U$ is one-to-one.
Let $z_1,z_2\in U$ such that $\tilde{\sigma}(z_1)=\tilde{\sigma}(z_2)$.
This means, by def\/inition, that
\begin{gather}
\label{eq:equal}
(\sigma(\rho(z_1)),-1,\rho(z_1))\cdot z_1=(\sigma(\rho(z_2)),-1,\rho(z_2))\cdot z_2.
\end{gather}
Since $\rho((\sigma(\rho(z_1)),-1,\rho(z_1))\cdot z_1)=\sigma(\rho(z_1))$ and
$\rho((\sigma(\rho(z_2)),-1,\rho(z_2))\cdot z_2)=\sigma(\rho(z_2))$, and since $\rho(z_1)$ and $\rho(z_2)$ are elements
of~$V$, it follows that $\rho(z_1)=\rho(z_2)$.
Moreover, if we multiply the equation~\eqref{eq:equal} on the left~by
\begin{gather*}
(\sigma(\rho(z_1)),-1,\rho(z_1))^{-1}=(\rho(z_1),1,\sigma(\rho(z_1))),
\end{gather*}
we obtain that
\begin{gather*}
z_1=(\rho(z_1),1,\sigma(\rho(z_1)))\cdot \bigl((\sigma(\rho(z_2)),-1,\rho(z_2))\cdot z_2\bigr)
\\
\phantom{z_1}
 =(\rho(z_1),0,\rho(z_2))\cdot z_2=(\rho(z_2),0,\rho(z_2))\cdot z_2
 =z_2.
\end{gather*}
Thus $\tilde{\sigma}|_U$ is one-to-one.
Observe that $\rho(U)=V$ and $\sigma(V)$ are both open.
Hence, $\mathcal{G}(\sigma(V),0,1$, $V) * U$ is open.
Since the action map~$\alpha$ is open by Lem\-ma~\ref{lem:etale},
\begin{gather*}
\tilde{\sigma}(U)= \alpha\big(\mathcal{G}(\sigma(V),0,1,V) * U\big)
\end{gather*}
is open.
Therefore $\tilde{\sigma}$ is a~local homeomorphism.
The equation $\rho\circ \tilde{\sigma}=\sigma\circ\rho$ is an easy computation.

We now prove the second part of the theorem, namely that $\Psi:\tilde{\mathcal{G}}\to \mathcal{G}\ast Z$ (see
formula~\eqref{eq:isomorphism}) is an isomorphism where $\tilde{\mathcal{G}} = \mathcal{G}(Z, \tilde{\sigma})$ is the
Renault--Deaconu groupoid associated to $\tilde{\sigma}$, that is,
\begin{gather*}
\tilde{\mathcal{G}}=\big\{ (t,m-n,s)\in Z\times \mathbb{Z}\times Z: \tilde{\sigma}^m(t)=\tilde{\sigma}^n(s) \big\}.
\end{gather*}
We f\/irst show that~$\Psi$ is well def\/ined.
Let $(t,m-n,s)\in \tilde{\mathcal{G}}$; then $\tilde{\sigma}^m(t)=\tilde{\sigma}^n(s)$.
An easy computation shows that $\tilde{\sigma}^m(t)=(\sigma^m(\rho(t)),-m,\rho(t))\cdot t$.
Thus, if $\tilde{\sigma}^m(t)=\tilde{\sigma}^n(s)$ then $\sigma^m(\rho(t))=\sigma^n(\rho(s))$ and, moreover,
$t=(\rho(t),m-n,\rho(s))\cdot s$.
Hence $(\rho(t),m-n,\rho(s))\in \mathcal{G}$ and~$\Psi$ is well def\/ined.

We next show that~$\Psi$ is a~groupoid morphism.
Let $(t_1,n_1,s_1)$ and $(t_2,n_2,s_2)$ be two composable elements in $\tilde{\mathcal{G}}$.
Then $s_1=t_2$ and
\begin{gather*}
(t_1,n_1,s_1)(t_2,n_2,s_2)=(t_1,n_1+n_2,s_2).
\end{gather*}
Therefore
\begin{gather}
\label{eq:1}
\Psi((t_1,n_1,s_1)(t_2,n_2,s_2))=((\rho(t_1),n_1+n_2,\rho(s_2)),s_2).
\end{gather}
Moreover, since $s_1=t_2=(\rho(t_2),n_2,\rho(s_2))\cdot s_2$, it follows that $((\rho(t_1),n_1,\rho(s_1)),s_1)$ and
$((\rho(t_2),n_2$, $\rho(s_2)),s_2)$ are composable and we have that
\begin{gather}
\Psi((t_1,n_1,s_1)) \Psi((t_2,n_2,s_2)) = ((\rho(t_1),n_1,\rho(s_1)),s_1)\cdot ((\rho(t_2),n_2,\rho(s_2)),s_2)
\nonumber
\\
\hphantom{\Psi((t_1,n_1,s_1)) \Psi((t_2,n_2,s_2))}{}
=((\rho(t_1),n_1+n_2,\rho(s_2)),s_2).
\label{eq:2}
\end{gather}
Equations~\eqref{eq:1} and~\eqref{eq:2} imply that~$\Psi$ is a~groupoid morphism.

Next we prove that~$\Psi$ is one-to-one and onto.
Let $(t_1,n_1,s_1)$ and $(t_2,n_2,s_2)\in \tilde{\mathcal{G}}$ such that $\Psi((t_1,n_1,s_1))=\Psi((t_2,n_2,s_2))$.
Thus
\begin{gather*}
\bigl((\rho(t_1),n_1,\rho(s_1)),s_1\bigr)=\bigl((\rho(t_2),n_2,\rho(s_2)),s_2\bigr).
\end{gather*}
Therefore $s_1=s_2$, $\rho(t_1)=\rho(t_2)$, and $n_1=n_2$.
Moreover
\begin{gather*}
t_1=(\rho(t_1),n_1,\rho(s_1))\cdot s_1 =(\rho(t_2),n_2,\rho(s_2))\cdot s_2 =t_2.
\end{gather*}
Hence~$\Psi$ is one-to-one.
To see that~$\Psi$ is onto, let $((x,n,y),s)\in \mathcal{G}\ast Z$, and recall that $y=\rho(s)$.
Let $t=(x,n,y)\cdot s$.
Then $\rho(t)=x$ and
\begin{gather*}
\Psi(t,n,s)=((\rho(t),n,\rho(s)),s)=((x,n,y),s).
\end{gather*}

To show that~$\Psi$ is continuous, let $U=\mathcal{G}(U_1, m, n, U_2) \ast V_2$ be an open set in $\mathcal{G}\ast Z$
such that $U_1$, $U_2$ are open subsets of~$X$, $\sigma^m(U_1)=\sigma^n(U_2)$ is open, $\sigma^m|_{U_1}:U_1\to
\sigma^m(U_1)$ and $\sigma^n|_{U_2}:U_2\to \sigma^n(U_2)$ are homeomorphisms, and $V_2=\rho^{-1}(U_2)$.
Then
\begin{gather*}
\Psi^{-1}(U)=\big(\rho^{-1}(U_1), m, n,\rho^{-1}(U_2)\big).
\end{gather*}
One can repeat the argument that we made when showing that $\tilde{\sigma}$ is a~local homeomorphism and prove that
$\tilde{\sigma}^m(\rho^{-1}(U_1))$ is open and $\tilde{\sigma}^m|_{\rho^{-1}(U_1)}:\rho^{-1}(U_1)\to
\tilde{\sigma}^m(\rho^{-1}(U_1))$ is a~homeomorphism, and, similarly, $\tilde{\sigma}^n(\rho^{-1}(U_2))$ is open and
$\tilde{\sigma}^n|_{\rho^{-1}(U_2)}:\rho^{-1}(U_2)\to \tilde{\sigma}^n(\rho^{-1}(U_2))$ is a~homeomorphism.
Notice that $\tilde{\sigma}^m(t)=\tilde{\sigma}^n(s)$ with $t\in \rho^{-1}(U_1)$ and $s\in \rho^{-1}(U_2)$ if and only
if $\sigma^m(\rho(t))=\sigma^n(\rho(s))$.
Since $\sigma^m(U_1)=\sigma^n(U_2)$ it follows that $\tilde{\sigma}^m(\rho^{-1}(U_1))=\tilde{\sigma}^n(\rho^{-1}(U_2))$.
Therefore $\Psi^{-1}(U)$ is an open cylinder in $\tilde{\mathcal{G}}$ and~$\Psi$ is continuous.
Finally, $\Psi^{-1}$ is continuous since it is a~composition of continuous maps.
\end{proof}

For the rest of this section we assume that both~$X$ and~$Z$ are locally compact spaces.
This assumption implies that both $\mathcal{G}$ and $\tilde{\mathcal{G}}$ are locally compact groupoids.

\begin{Corollary}
\label{cor:amenable}
Assume that~$\sigma$ is a~local homeomorphism on a~locally compact space~$X$ and assume that the Renault--Deaconu
groupoid $\mathcal{G}=\mathcal{G}(X,\sigma)$ acts on the locally compact space~$Z$.
The groupoid $\tilde{\mathcal{G}} \cong \mathcal{G}\ast Z$ is amenable and the full and reduced $C^*$-algebra of
$\mathcal{G}\ast Z$ coincide.
Moreover, $C^*(\mathcal{G}\ast Z)$ is nuclear.
\end{Corollary}
\begin{proof}
This is an immediate consequence of~\cite[Proposition 2.9]{Ren_CLA00}.
\end{proof}

As another application of Theorem~\ref{thm:lift}, we can identify $C^*(\mathcal{G}\ast Z)$ as a~Cuntz--Pimsner algebra~\cite{mPi_FIC95}.
The closure of $C_c(Z)$ under a~suitable norm may be viewed as a~$C^*$-correspondence $\mathcal{Z}$ over the
$C^*$-algebra $C_0(Z)$.
The inner product is def\/ined for $\xi,\eta \in C_c(Z)$~by
\begin{gather*}
\langle\xi,\eta\rangle(z)= \sum\limits_{\tilde{\sigma} (y)=z}\overline{\xi(y)}\eta(y);
\end{gather*}
the left and right actions of $C_0(Z)$ are given by $(a\cdot\xi\cdot b)(z)=a(z)\xi(z)b(\tilde{\sigma}(z))$ where $a, b
\in C_0(Z)$ and $\xi\in C_c(Z)$.
The following corollary follows from~\cite[Theorem 7]{De_Ku_Mu_JOT01}.
\begin{Corollary}
Under the hypothesis of Corollary~{\rm \ref{cor:amenable}}, the groupoid $C^*$-algebra $C^*(\mathcal{G}\ast Z)$ is isomorphic
to the Cuntz--Pimsner $C^*$-algebra $\mathcal{O}_{\mathcal{Z}}$.
\end{Corollary}

\section{Groupoid actions on fractafolds}
\label{sec:groupoidaction}

In this section we build a~fractafold bundle~$L$ such that the Renault--Deaconu $\mathcal{G}$ groupoid associated to
a~one-sided shift map~$\sigma$ acts naturally on~$L$.
The action groupoid $\mathcal{G}*L$ encapsulates the symmetries of the fractafold bundle.
Using the results of the previous section,~$\sigma$ extends to a~local homeomorphism $\tilde{\sigma}$ on~$L$.
We prove that $\tilde{\sigma}$ is essentially free and, hence, $\mathcal{G}*L$ is topologically principal.
Moreover, we show that the action groupoid $\mathcal{G}*L$ contains a~dense orbit.
If~$L$ is locally compact, it follows that $C^*(\mathcal{G}*L)$ is primitive.
If the iterated function system def\/ining~$L$ satisf\/ies the open set condition, then we build a~$\mathcal{G}*L$-invariant
measure on~$L$.
Thus if~$L$ is, in addition, locally compact, then $C^*(\mathcal{G}*L)$ has a~densely def\/ined lower semicontinuous
trace.

Let $(T,d)$ be a~complete metric space and let $(F_1,\dots, F_{N})$ be an \emph{iterated function system} on~$T$~\cite{Bar_FE,Edg_MTFG2, Hut_81}.
That is, each $F_i$ is a~strict contraction on~$T$.
We assume that the iterated function system is \emph{non-degenerate},
meaning that there are constants $0<r_i\le R_i<1$, $i=1,\dots, N$, such that
\begin{gather}
\label{eq:nondegifs}
r_id(t,t^\prime)\le d(F_i(t),F_i(t^\prime))\le R_id(t,t^\prime)
\qquad
\text{for all}
\quad
t,
t^\prime\in T,
\qquad
i=1,\dots, N.
\end{gather}
We further assume map $F_i$ is surjective for all $i=1,\dots, N$; so the $F_i$ are homeomorphisms.
If we can chose $r_i=R_i$ in~\eqref{eq:nondegifs}, then $F_i$ is a~similarity (or a~similitude).
For an iterated function system there
is a~unique compact \emph{invariant set}~$K$ \cite[Theorem~3.1.3]{Hut_81} such that
\begin{gather*}
K=F_1(K)\cup \dots\cup F_{N}(K).
\end{gather*}
If $F_i(K)\bigcap F_j(K)=\varnothing$ for all $i\ne j$, then the iterated function system is called \emph{totally
disconnected}.
In this case the invariant set~$K$ is a~totally disconnected set \cite[Theorem~8.2.1]{Bar_FE}.
If $N\ge 3$, using~\cite[Remark~3.1.9]{Hut_81} one can f\/ind examples of iterated function systems that are not totally
disconnected but have a~totally disconnected invariant set.

Let $W=\{1,\dots,N\}$ and def\/ine $W^*=\bigcup_{n\in\mathbb{N}}W^n$ to be the set of f\/inite words over the alphabet~$W$,
and $X=W^\infty$ to be the set of inf\/inite words (sequences) with elements in~$W$.
If $\omega\in W^n$ we say that the length of~$\omega$, denoted by $\vert\omega\vert$, is~$n$.
For $\omega\in W^n$ we write $F_\omega(x):=F_{\omega_{n}}\circ \dots \circ F_{\omega_1}(x)$,
$F_{\omega}^{-1}(x):=F_{\omega_{1}}^{-1}\circ \dots\circ F_{\omega_{n}}^{-1}(x)$, $r_\omega:=r_{\omega_1}\cdots
r_{\omega_{n}}$ and similarly for $R_\omega$.

If $\omega\in W^n$ or $x\in X$ and if $k\in\mathbb{N}$ such that $k\le n$ (if applicable), then we write
\begin{gather*}
\omega(k):=\omega_1\cdots \omega_{k}
\qquad
\text{and}
\qquad
x(k):=x_1\cdots x_{k}.
\end{gather*}

Given a~f\/inite sequence $\omega\in W^n$ or an inf\/inite sequence $x\in X$ set
\begin{gather*}
L_n(\omega)=F_{\omega(n)}^{-1}(K)
\qquad
\text{and}
\qquad
L_n(x):=L_n(x(n)).
\end{gather*}
Then, if $x\in X$, $L_n(x)\subset L_{n+1}(x)$ and we def\/ine $L(x):= \bigcup_{n\in\mathbb{N}}L_n(x)$ endowed with the
inductive limit topology.
We refer to $L(x)$ as the \emph{infinite blow-up} of~$K$ at~$x$.
Our def\/inition is not quite the same as the one used by Stricharz since he uses the relative topology
(see~\cite{Str_CJM98} and \cite[Chapter~5.4]{Str_Prin06}); but the two seem to agree in some cases, for example, when the
blowups are zero dimensional.
Two blow-ups $L(x)$ and $L(y)$ are homeomorphic if the inf\/inite sequences~$x$ and~$y$ dif\/fer in a~f\/inite number of
indices.
For example, if $y = x_2x_3\cdots$, $F_{x_1}$ extends to a~homeomorphism from $L(x)$ to $L(y)$.

If~$U$ is a~subset of~$K$ and $\omega\in W^n$ is a~f\/inite sequence or $x\in X$ is an inf\/inite word, then we write
$L_n^\omega(U)$ for $F_{\omega(n)}^{-1}(U)$ and $L_n^x(U)$ for $L_n^{x(n)}(U)$.

If $\omega\in W^*$, then the clopen cylinder $Z(\omega)\subset X=W^\infty$ is def\/ined via
\begin{gather*}
Z(\omega)=\{x\in X:x_i=\omega_i, \, i=1,\dots,\vert\omega\vert\}.
\end{gather*}
The collection $\{Z(\omega)\}_{\omega\in W^*}$ forms a~basis of a~topology on~$X$, and, endowed with this topology,
$X$ is a~totally disconnected compact space.
Moreover, the shift map $\sigma:X\to X$ def\/ined by $\sigma(x_1x_2\cdots)=(x_2x_3\cdots)$ is a~local homeomorphism on~$X$.
We write $\mathcal{G}$ for the Renault--Deaconu groupoid associated to~$\sigma$ as in Section~\ref{sec:Groupoidactions}.
Recall that the unit space $\mathcal{G}^{0}$ is homeomorphic to~$X$.

Next we build a~fractafold bundle on which the groupoid $\mathcal{G}$ acts.
For $n\ge 0$ def\/ine
\begin{gather*}
L_n=\bigsqcup_{\omega\in W^n}Z(\omega)\times L_n(\omega)\subset X\times Y.
\end{gather*}
Then each $L_n$ is a~compact space and $L_n\subset L_{n+1}$.
Observe that that $L_0=X\times K$.
We def\/ine the \emph{fractafold bundle}~$L$ to be the increasing union of $L_n$, $L=\bigcup_{n\ge 0}L_n$, endowed with
the inductive limit topology.
That is, a~set $U\subset L$ is open if and only if $U\cap L_n$ is open for all $n\ge 0$.
We will show below that the base space is~$X$.
The following characterization of open sets is used later in this section.
\begin{Lemma}
Let $x\in X$ and $n\ge 0$.
Then $U\subset L_n(x)$ is open if and only if there is an open set $V\subset K$ such that $U=L_n^x(V)$.
\end{Lemma}
\begin{proof}
Recall that the maps $F_i$, $i=1,\dots, N$, are homeomorphisms.
Therefore the maps~$F_{\omega}$ and~$F_{\omega}^{-1}$ are homeomorphisms for all $\omega\in W^*$.
The conclusion of the lemma follows immedia\-tely.~~
\end{proof}

Recall that, in general, the inductive limit of an increasing sequence of Hausdorf\/f spaces might fail to be Hausdorf\/f.
However, as we prove in the following lemma, the bundle~$L$ is Hausdorf\/f.
\begin{Lemma}
The bundle~$L$ endowed with the inductive limit topology as above is a~Hausdorff space.
\end{Lemma}
\begin{proof}
Let $\iota:L\to X\times T$ be the inclusion map.
If~$U$ is an open set in $X\times T$, then $\iota^{-1}(U)\cap L_n=U\cap L_n$ is open in $L_n$ for all $n\ge 1$.
Therefore~$\iota$ is a~continuous one-to-one map.
Since $X\times T$ is Hausdorf\/f, it follows that~$L$ is Hausdorf\/f as well.
\end{proof}

Before we def\/ine the action of the Renault--Deaconu groupoid on~$L$ we show that the natural projection from~$L$
into~$X$ is an open map.
\begin{Lemma}
The map $\pi:L\to X$, $\pi(x,t)=x$ is a~continuous open map.
Moreover, $\pi^{-1}(\{x\})$ is homeomorphic to the fractafold $L(x)$.
\end{Lemma}
\begin{proof}
Since the restriction of~$\pi$ to $L_n$ is continuous for each~$n$ and~$L$ is endowed with the inductive limit
topology,~$\pi$ is also continuous.
Let~$U$ be an open set in~$L$.
Therefore $U_n:=U\cap L_n$ is open for all $n\ge 0$.
Hence $U_n\cap(Z(\omega)\times L_n(\omega))$ is open for all $\omega\in W^n$.
For $\omega\in W^n$, the map $\pi|_{Z(\omega)\times L_n(\omega)}:Z(\omega)\times L_n(\omega)\to Z(\omega)$ is open
because it is just the projection onto the f\/irst coordinate.
Since $Z(\omega)$ is an open subset of~$X$ it follows that $\pi(U_n\cap(Z(\omega)\times L_n(\omega)))$
is open in~$X$ for all $\omega\in W^n$ and $n\ge 0$.
Therefore, $\pi(U_n)$ is open in~$X$ for all $n\ge 0$, and hence $\pi(U) = \bigcup_n \pi(U_n)$ is open.
Thus~$\pi$ is an open map.
The f\/inal assertion is obvious.
\end{proof}

The next results shows that the Renault--Deaconu groupoid $\mathcal{G}$ associated to the shift map on~$X$ acts on the
left on the space~$L$.
Note that if $\gamma=(x,m-n,y)\in \mathcal{G}$ and $(z,t)\in L$, then $s(\gamma)=\pi(z,t)$ if and only if $y=z$.
Let
\begin{gather}
\label{eq:actiongroupoid}
\mathcal{G}\ast L=\{((x,m-n,y),(y,t)):(x,m-n,y)\in\mathcal{G},(y,t)\in L\}.
\end{gather}

\begin{Theorem}
\label{thm:action}
With the notation as above, the Renault--Deaconu groupoid $\mathcal{G}$ associated to the shift map~$\sigma$ on~$X$ acts
on the fractafold bundle~$L$ via the map $((x,m-n,y),(y,t))\mapsto (x,m-n,y) \cdot (y,t)$ defined on $\mathcal{G}\ast
L$, where
\begin{gather}
\label{eq:action}
(x,m-n,y) \cdot (y,t) = \big(x,F_{x(m)}^{-1}\circ F_{y(n)}(t)\big).
\end{gather}
Moreover, the action map is open.
\end{Theorem}
\begin{proof}
We need to prove that the above map is well def\/ined; that is, we need to show that the range of the map is~$L$.
Let $(x,m-n,y)\in\mathcal{G}$ and $(y,t)\in L$.
Let $k\ge 0$ such that $(y,t)\in L_k$.
Then there is $\omega\in W^k$ such that $y\in Z(\omega)$ and $t\in L_k(\omega)$.
Therefore $y_i=\omega_i$ for $i=1,\dots, k$.
Notice that it suf\/f\/ices to assume that $n\ge k$.
Indeed, if $m^\prime-n^\prime =m-n$ and $\sigma^{m^\prime}(x)=\sigma^{n^\prime}(y)$ with $m^\prime>m$ and $n^\prime>n$,
then $x_{m+i} = y_{n+i}$ for all $i=1,\dots,m^\prime-m$ and so $F_{x(m)}^{-1}\circ F_{y(n)}(t) =
F_{x(m^\prime)}^{-1}\circ F_{y(n^\prime)}(t)$.
Then
\begin{gather*}
F_{y(n)}(t)\in K
\end{gather*}
and, thus, $F_{x(m)}^{-1}\circ F_{y(n)}(t)\in L_m(\alpha)$, where $\alpha=x(m)$.
Hence{\samepage
\begin{gather*}
\big(x,F_{x(m)}^{-1}\circ F_{y(n)}(t)\big)\in Z(\alpha)\times L_m(\alpha)\subset L.
\end{gather*}
So the action is well def\/ined.}

The left action is continuous since for f\/ixed~$m$ and~$n$ in $\mathbb{N}$
and f\/ixed words~$\alpha\in W^m$ and \mbox{$\beta\in W^n$}, we have that the map $F_\alpha^{-1}\circ F_\beta$ is continuous.
Thus the Renault--Deaconu groupoid~$\mathcal{G}$ acts on the left on~$L$.
The last part of the theorem is an immediate consequence of Lem\-ma~\ref{lem:etale}.~~~
\end{proof}

\begin{Corollary}
\label{cor:lift}
The map $\tilde{\sigma}:L\to L$ defined by $\tilde{\sigma}(x,t)=(\sigma(x),F_{x_1}(t))$ is a~local homeomorphism on~$L$
such that $\pi\circ\tilde{\sigma}=\sigma\circ\pi$.
\end{Corollary}
\begin{proof}
Theorems~\ref{thm:lift} and~\ref{thm:action} imply that the shift map~$\sigma$ on~$X$, which is a~local homeomorphism,
lifts to a~local homeomorphism $\tilde{\sigma}$ on~$L$ such that $\tilde{\sigma}(x,t)=(\sigma(x),-1,x)\cdot (x,t)$.
By equation~\eqref{eq:action}, $\tilde{\sigma}(x,t)=(\sigma(x),F_{x_1}(t))$.
\end{proof}

\begin{Corollary}
Let $\tilde{\sigma}$ be the local homeomorphism on~$L$ provided by Corollary~{\rm \ref{cor:lift}}.
Then the left action groupoid $\mathcal{G}\ast L$ is homeomorphic to the Renault--Deaconu groupoid $\tilde{\mathcal{G}}$
associated to $\tilde{\sigma}$.
\end{Corollary}
\begin{proof}
This is an immediate consequence of the second part of Theorem~\ref{thm:lift}.
\end{proof}

Recall from~\cite[Def\/inition on p.~1781]{Dea_TAMS95} that a~local homeomorphism~$\tau$ on a~topological space~$Z$ is
\emph{essentially free} if
\begin{gather*}
\{ z\in Z
\mid
\forall k,l\ge 0,\tau^k(z)=\tau^l(z)\Rightarrow k=l\}
\end{gather*}
is dense in~$Z$.
If~$\tau$ is essentially free then the Renault--Deaconu groupoid $\mathcal{G}(Z,\tau)$ is topologically principal
because, for $z\in Z$, the isotropy group $\mathcal{G}(z)$ is nontrivial if and only if
there are $k,l\ge 0$ with $k\ne l$ such that $\tau^k(z) = \tau^l(z)$ \cite[Example~1.2c]{ADC_BSMF97}.
It is easy to see that the shift map~$\sigma$ on~$X$ is essentially free \cite[Example~2]{Dea_TAMS95}.

\begin{Proposition}
\label{prop:essen_free}
The local homeomorphism $\tilde{\sigma}$ on~$L$ defined in Corollary~{\rm \ref{cor:lift}} is essentially free.
Hence, the Renault--Deaconu groupoid $\tilde{\mathcal{G}}$ associated to $\tilde{\sigma}$ is topologically principal.
\end{Proposition}
\begin{proof}
Let~$U$ be a~nonempty open subset of~$L$.
We need to f\/ind $(x,t)\in U$ such that for all $k,l\ge 0$ if $\tilde{\sigma}^k(x,t)=\tilde{\sigma}^l(x,t)$ then $k=l$.
Corollary~\ref{cor:lift} implies that if $\tilde{\sigma}^k(x,t)=\tilde{\sigma}^l(x,t)$ then $\sigma^k(x)=\sigma^l(x)$.
Since $\pi (U)$ is open in~$X$ there is $x\in \pi(U)$ such that for all $k,l\ge 0$ if $\sigma^k(x)=\sigma^l(x)$ then
$k=l$.
It follows that if we pick $t\in L_x$ such that $(x,t)\in U$ then $(x,t)$ satisf\/ies the desired property.
Hence $\tilde{\sigma}$ is essentially free and $\tilde{\mathcal{G}}$ is topologically principal.
\end{proof}

\begin{Proposition}
\label{prop:dense_orbit}
Let $\mathcal{G}$ be the Renault--Deaconu groupoid associated to the shift map~$\sigma$ on~$X$ and let~$L$ be the
fractafold bundle associated to a~non-de\-ge\-ne\-rate iterated function system $(F_1,\dots, F_N)$ on a~complete metric
space~$Y$.
Let $\mathcal{G}\ast L$ be the left action groupoid defined via~\eqref{eq:actiongroupoid} and~\eqref{eq:action}.
Let $x\in X$ be an infinite word obtained by concatenating all the finite words in $W^*$. Then for all $t \in K \subset
L_x$, the orbit of $(x,t)\in L\simeq(\mathcal{G}\ast L)^0$ is dense.
\end{Proposition}
\begin{proof}
Recall (see, for example,~\cite[Theorem 4.2.1]{Bar_FE}) that every point~$v$ in~$K$ has at least one address $y\in X$,
that is, $\{v\}=\bigcap_{n\ge 0}F_{y_{1}}\cdots F_{y_{n}}(K)$.
Notice that the diameter of $F_{y_1}\cdots F_{y_{n}}(K)$ is less than $R_{y(n)}\cdot\operatorname{diam}K$ and recall that
each $R_i$ is strictly smaller than~$1$.
We claim that the sequence $\{F_{x_{n}}\cdots F_{x_{1}}(t)\}_{n\in\mathbb{N}}$ is dense in~$K$ for all $t\in K$.
To prove the claim, let $v\in K$ and $\varepsilon>0$.
There is $k\in\mathbb{N}$ and $\omega\in W^k$ such that $d(v,F_{\omega_1}\cdots F_{\omega_{k}}(u))<\varepsilon$ for all $u\in K$.
By hypothesis,~$x$ contains the word $\omega_{k}\cdots\omega_1$.
That is, there is $l\ge 1$ such that $x_l=\omega_{k}, \dots$, $x_{l+k}=\omega_1$.
Then
\begin{gather*}
d(v,F_{x_{l+k}}\cdots F_{x_l}(F_{x_{l-1}}\cdots F_{x_1}(t))) =d(v,F_{\omega_{1}}\cdots F_{\omega_{k}}(F_{x_{l-1}}\cdots
F_{x_1}(t)))<\varepsilon
\end{gather*}
for all $t\in K$.
The claim follows.

We prove that the orbit of $(x,t)\in L_0$ is dense, where~$t$ is an arbitrary point in~$K$.
Let $(y,v)\in L$ and let~$U$ be a~neighborhood of $(y,v)$ in~$L$.
Let $m\in \mathbb{N}$ and $\alpha\in W^m$ such that $(y,v)\in U_m:=U\cap(Z(\alpha)\times L_{m}(\alpha))$.
Then there is~$V$ open in~$K$ such that $(y,v)\in Z(\alpha)\times L_{m}^\alpha(V)$.
We need to f\/ind $\gamma\in \mathcal{G}$ such that $\gamma\cdot (x,t)\in U$.
Let $n\in\mathbb{N}$ be such that $F_{x(n)}(t)\in V$.
Then $F_{\alpha(m)}^{-1}\circ F_{x(n)}(t)\in L_m^\alpha(V)$.
Def\/ine $y\in X$ such that $y_i=\alpha_i$, $i=1,\dots,m$ and $y_{m+i}=x_{n+i}$ for all $i\ge 1$.
Then $(y,m-n,x)\in \mathcal{G}$ and $(y,m-n,x)\cdot (x,t)\in U$.
\end{proof}

In general~$L$ may fail to be locally compact (see Section~\ref{sec:examples} for specif\/ic examples of when this
property fails).
However, if $F_i(K)$ is open in~$K$ for all $i=1,\dots, N$, then~$L$ is locally compact.
In this case~$K$ is a~totally disconnected set.
One can easily check that this condition is satisf\/ied if the iterated function system is totally disconnected.

Suppose that~$L$ is locally compact, then as observed at the end of Section~\ref{sec:Groupoidactions},
$C^*(\mathcal{G}*L)$ is nuclear and a~Cuntz--Pimsner algebra
(where the correspondence is def\/ined over $C_0(L)$).
Proposition~\ref{prop:dense_orbit} allows us to conclude a~bit more about $C^*(\mathcal{G} *L)$.

A $C^*$-algebra is def\/ined to be \emph{primitive} (see~\cite[\S\,3.13.7]{gkp_79}) if it has a~faithful irreducible
representation.

\begin{Proposition}
Assume that $F_i(K)$ is open in~$K$ for all $i=1,\dots,N$.
With notation as in Proposition~{\rm \ref{prop:dense_orbit}}, $C^*(\mathcal{G}*L)$ is primitive.
\end{Proposition}
\begin{proof}
Let $z = (x, t) \in L\simeq(\mathcal{G}\ast L)^0$ be a~point with dense orbit as guaranteed~by
Proposition~\ref{prop:dense_orbit}.
Note that point evaluation at~$z$ def\/ines a~pure state on~$C_0(L)$, the canonical masa (maximal abelian subalgebra) in~$C^*(\mathcal{G}*L)$.
Since it extends to a~pure state on $C^*(\mathcal{G}*L)$, the GNS construction provides an irreducible representation
$\pi_z$ of $C^*(\mathcal{G}*L)$ and a~unit vector $\xi_z$ in the associated Hilbert space $H_{\pi_z}$ such that $\langle
\pi_z(f)\xi_z, \xi_z \rangle = f(z)$ for all $f \in C_0(L)$.
Since $\ker\pi_z \cap C_0(L)$ is supported on an open set which does not contain any points in the orbit of~$z$,
$\ker\pi_z \cap C_0(L) = \{0 \}$.
By Corollary~\ref{cor:amenable} and Proposition~\ref{prop:essen_free}, $\mathcal{G}*L = \tilde{\mathcal{G}}$ is both
amenable and topologically principal.
Hence, every nonzero ideal in $C^*(\mathcal{G}*L)$ must have a~nontrivial intersection with $C_0(L)$ by~\cite[Theorem~4.4]{MR2745642}
(see also~\cite[Proposition~5.5]{BCFS-ppt12}).
Hence, $\ker\pi_z = \{0 \}$ and thus $C^*(\mathcal{G}*L)$ has a~faithful irreducible representation.
\end{proof}

If~$G$ is an \'etale groupoid, a~$G$-invariant measure~$\mu$ on $G^0$ is a~measure such that for any open~$G$-set~$U$,
$\mu(r(U))=\mu(s(U))$ (it suf\/f\/ices to prove this for basis of open~$G$-sets).
We show below that there is a~$\mathcal{G}*L$-invariant measure $\mu_\infty$ on~$L$ provided that the iterated function
system satisf\/ies one additional hypothesis.
We will assume that the iterated function system $(F_1,\dots,F_N)$ satisf\/ies the \emph{open set condition}
\cite[Def\/inition~5.2.1]{Hut_81}: there exists a~non-empty open set~$O$ such that $\cup_i F_i(O)\subset O$ and
$F_i(O)\cap F_j(O)=\varnothing$ if $i\ne j$.
A~totally disconnected iterated function system satisf\/ies the open set condition with $O=K$.

We will need to invoke a~measure theoretic extension theorem (see~\cite[Theorem 6.2]{MW13})
to extend a~measure on a~semialgebra of subsets of~$L$ to the~$\sigma$-algebra of Borel sets on~$L$.
Recall that a~collection $\mathcal{C}$ of subsets of a~set~$\Omega$ is called a~\emph{semialgebra} if it is closed under
f\/inite intersections and if the complement of $B \in \mathcal{C}$ is expressible as a~f\/inite disjoint union of elements
of $\mathcal{C}$.
For each $n\in\mathbb{N}$, let $\mathcal{C}_n$ be the collection of Borel subsets of $L_n$.
Then $\mathcal{C} = \{L_n^c: n \ge 1 \} \cup \big(\bigcup_{n \in \mathbb{N}} \mathcal{C}_n\big)$ is a~semialgebra.

By~\cite[Theorem 4.4.1]{Hut_81}, there is a~unique Borel probability measure~$\mu$ on~$K$ such that
\begin{gather*}
\mu(A)=\frac{1}{N}\sum\limits_{i=1}^N \mu\big(F_i^{-1}(A)\big)
\end{gather*}
for all Borel subsets of~$K$.
Then, for each $n\in \mathbb{N}$ and $\omega\in W^n$ we can def\/ine a~measure $\mu_\omega$ on $L_n(\omega)$ via
$\mu_\omega(A)=N^n\mu(F_\omega(A))$.
We let~$\nu$ be the product measure on~$X$ generated by the weights $\{1/N,\dots,1/N\}$ on the set $\{1,\dots,N\}$.
Therefore, if $\omega\in W^n$ then $\nu(Z(\omega))=(1/N)^n$.
Then one can def\/ine a~Borel measure $\mu_n$ on $L_n$ such that $\mu_n(Z(\omega)\times A)=\mu(F_\omega(A))$, and, more
generally, $\mu_n(U\times A)=\nu(\sigma^n(U))\cdot \mu(F_\omega(A))$ if $U\times A\subset Z(\omega)\times L_n(\omega)$.
Notice that if $m<n$ and~$B$ is a~Borel subset of $L_m$ then~$B$ is a~Borel subset of $L_n$.
However, more is true if the iterated function system satisf\/ies the open set condition.
\begin{Lemma}
\label{lem:osc}
Assume that the iterated functions system $(F_1,\dots,F_N)$ on~$T$ satisfies the open set condition.
Then
\begin{enumerate}\itemsep=0pt
\item[$1.$] The measure~$\mu$ is the Hausdorff measure on~$K$ and $\mu(A) = N\mu(F_i(A))$ for all Borel sets $A \subset K$.
\item[$2.$] If $m<n$ and~$B$ is a~Borel subset of $L_m$ then $\mu_n(B)=\mu_m(B)$.
\end{enumerate}
\end{Lemma}
\begin{proof}
The f\/irst part is an immediate consequence of~\cite[Theorem 5.3.1]{Hut_81}.
The second part follows from the f\/irst part.
\end{proof}

\begin{Proposition}
Assume that the iterated function system $(F_1,\dots,F_N)$ satisfies the open set condition.
Then there is a~unique $\tilde{\mathcal{G}}$-invariant Borel measure $\mu_\infty$ on~$L$ such that $\mu_\infty(L_n^c) =
\infty$ and
\begin{gather*}
\mu_\infty(U \times A) = \nu(\sigma^n(U))\cdot\mu(F_\omega(A))
\end{gather*}
where $\omega \in W^n$ and $U\times A \subset Z(\omega)\times L_n(\omega)$ is Borel.

\end{Proposition}
\begin{proof}
As noted above $\mathcal{C} = \{L_n^c: n \ge 1 \} \cup \big(\bigcup_{n \in \mathbb{N}} \mathcal{C}_n\big)$ is
a~semialgebra.
By~\cite[Theorem 6.2]{MW13}, $\mu_\infty$ extends uniquely to a~measure on the~$\sigma$-algebra generated~by
$\mathcal{C}$ if the following conditions hold
\begin{itemize}\itemsep=0pt
\item[i.] If $\varnothing \in \mathcal{C}$, $\mu_\infty(\varnothing) = 0$.
\item[ii.] If $C = C_1 \cup \dots \cup C_n$ where~$C$, $C_1, \dots, C_n \in \mathcal{C}$ and $C_i \cap C_j =
\varnothing$ for $i \ne j$, then $\mu_\infty(C) = \mu_\infty(C_1) + \dots + \mu_\infty(C_n)$.
\item[iii.] If $C \subset \bigcup_n C_n$ where~$C$, $C_1, C_2, \ldots \in \mathcal{C}$, then $\mu_\infty(C) \le
\sum\limits_n \mu_\infty(C_n)$.
\item[iv.] There are $C_1, C_2, \ldots \in \mathcal{C}$ such that $\mu_\infty(C_n) < \infty$ and $L = \bigcup_n C_n$.
\end{itemize}
Conditions (i) and (ii) are easy to check.
The only non-trivial case for condition (iii) is when $C=L_k^c$ for some $k\in \mathbb{N}$, and for each $n\in
\mathbb{N}$ there is $k_n\in \mathbb{N}$ such that $C_n\in \mathcal{C}_{k_n}$.
Then $\mu_\infty(C)=\infty$ and $\lim_{k\to \infty}\mu_\infty(L_k\bigcap C)=\infty$.
Let $k\in \mathbb{N}$ be f\/ixed.
Then the sets $L_k\cap C_n$ are Borel subsets of $L_k$ that cover $L_k\cap C$.
Hence, using Lemma~\ref{lem:osc}, we have that
\begin{gather*}
\mu_k(L_k\cap C)\le \sum\limits_j \mu_k(L_k\cap C_j)\le \sum\limits_n \mu_\infty (C_n).
\end{gather*}
Since the left hand side of the inequality goes to~$\infty$, it follows that $\sum\limits_n
\mu_\infty(C_n)=\infty=\mu_\infty(C)$.

To see that condition (iv) holds we f\/irst note that for each $\alpha \in W^n$, $C_\alpha = Z(\alpha) \times L_n(\alpha)
\in \mathcal{C}$, $\mu(C_\alpha) = \mu(K)=1 $ and $L = \bigcup_{\alpha \in W^*} C_\alpha$.
It is straightforward to show that the~$\sigma$-algebra generated by $\mathcal{C}$ is the Borels.

To prove that $\mu_\infty$ is $\tilde{\mathcal{G}}$-invariant we use the following fact that is probably known to
specialists: if~$\sigma$ is a~local homeomorphism on a~Hausdorf\/f topological space~$X$, then $\mathcal{G}(U,m,n,V)$ is
a~$\mathcal{G}$-section (where $\mathcal{G}=\mathcal{G}(X,\sigma)$ is the associated Renault--Deaconu groupoid), if and
only if~$\sigma^m|_U$ is a~homeomorphism onto $\sigma^m(U)$ and $\sigma^n|_V$ is a~homeomorphism onto $\sigma^n(V)$.
Moreover, $r(\mathcal{G}(U,m,n,V))=U$ and $s(\mathcal{G}(U,m,n,V))=V$.

Let~$V$ be open in~$L$ such that $\tilde{\sigma}|_V$ is a~homeomorphism onto $\tilde{\sigma}(V)$.
It follows from the above remark that in order to show that $\mu_\infty$ is $\tilde{\mathcal{G}}$-invariant it is enough
to prove that $\mu_\infty(V)=\mu_\infty(\tilde{\sigma}(V))$.
We have that
\begin{gather*}
\mu_\infty(V)=\mu(V\cap L_0)+\sum\limits_{n\ge 0}\mu_\infty(V\cap(L_{n+1}\setminus L_n))
\\
\phantom{\mu_\infty(V)}
=\mu(V\cap L_0)+\sum\limits_{n\ge 0}\sum\limits_{\omega\in W^{n+1}}\mu(V\cap(Z(\omega)\times (L_{n+1}(\omega)\setminus
L_n(\omega(n))))).
\end{gather*}
Hence, it suf\/f\/ices to prove that if $\omega\in W^n$ and $U\times A\subset Z(\omega)\times L_n(\omega)$ is open such that
$\tilde{\sigma}|_{U\times A}$ is a~homeomorphism onto $\tilde{\sigma}(U\times A)$, then $\mu(\tilde{\sigma}(U\times
A))=\mu(U\times A)$.
By the def\/inition of $\tilde{\sigma}$ we have that $\tilde{\sigma}(U\times A)=\sigma(U)\times F_{\omega_1}(A)$.
Hence, using Lemma~\ref{lem:osc},
\begin{gather*}
\mu_\infty(\tilde{\sigma}(U\times A))=\nu\big(\sigma^{n-1}(\sigma(U))\big)\cdot \mu(F_{\omega_n\cdots
\omega_2}(F_{\omega_1}(A)))=\mu_\infty(U\times A).
\end{gather*}
Thus $\mu_\infty$ is $\tilde{\mathcal{G}}$-invariant.
\end{proof}

\begin{Corollary}
Suppose that the iterated function system $(F_1,\dots,F_N)$ satisfies the open set condition and that $F_i(K)$ is open
in~$K$ for all $i=1,\dots, N$.
With notation as above, the map
\begin{gather*}
\tau(f)=\int_L fd\mu_\infty
\qquad
\text{for}
\quad
f\in C_c(\tilde{\mathcal{G}})\cap C^*(\tilde{\mathcal{G}})^+,
\end{gather*}
extends to a~densely defined lower semi-continuous trace on $C^*(\tilde{\mathcal{G}})\simeq C^*(\mathcal{G}*L)$.
\end{Corollary}

The proof of the corollary follows from the following lemma which is known to specialists.
We were unable, however, to f\/ind a~specif\/ic reference in the literature and we include a~short proof for completeness.

\begin{Lemma}
Suppose that~$G$ is an \'{e}tale locally compact groupoid and that~$\mu$ is a~$($Radon$)$~$G$-invariant measure on $G^0$.
Then the map
\begin{gather*}
\tau(f)=\int_{G^0}fd\mu
\qquad
\text{for all}
\quad
f\in C_c(G)\cap C^*(G)^+.
\end{gather*}
extends to a~densely defined lower semi-continuous trace on $C^*(G)$.
\end{Lemma}
\begin{proof}
Recall that for an \'{e}tale locally compact groupoid~$G$ the restriction map from~$C_c(G)$ to~$C_c(G^0)$ extends to
a~continuous expectation from $C^*(G)$ to $C^*(G^0)$ \cite[Proposition~II.1.15 and the remark following]{Ren_LNM793}.
It is well known that a~measure on a~locally compact space~$X$ induces a~lower semi-continuous weight on $C_0(X)$ (see,
for example,~\cite[Example~II.6.7.2(v)]{Bla_OA06}).
Hence we obtain a~densely def\/ined lower semi-continuous weight~$\tau$ on $C^*(G)$ by composition.
Finally, since~$\mu$ is invariant~\cite[Proposition~II.5.4]{Ren_LNM793} implies that~$\tau$ is a~trace.
\end{proof}

\section{Examples}
\label{sec:examples}

In this section we provide detailed descriptions of the fractafold bundle def\/ined in Section~\ref{sec:groupoidaction}
for some specif\/ic iterated function systems.
We point out some cases when the bundle~$L$ is not locally compact.
We show in the second example that the action groupoid $\mathcal{G}*L$ is not, in general, minimal.

\begin{Example}
\label{ex:unit}
Let $F_{0},F_{1}:\mathbb{R}\to\mathbb{R}$ be the maps $F_{0}(x)=\frac{1}{2}x$ and $F_{1}(x)=\frac{1}{2}x+\frac{1}{2}$.
Then $\{F_{0},F_{1}\}$ is an iterated function system whose invariant set is $[0,1]$.
In this example, $W=\{0,1\}$, $X=W^{\infty}=\{0,1\}^\infty$ and the fractafolds def\/ined in
Section~\ref{sec:groupoidaction} have an easy description.
Notice that, unlike in the previous section, we index our maps using $0$ and $1$; the index of an element in~$X$ starts at $0$ as well.
This makes the formulas that we describe next more tractable.
For $x\in X$ def\/ine $a_{n}(x)=F_{x(n)}^{-1}(0)$ and $b_{n}(x)=F_{x(n)}^{-1}(1)$ for $n\ge 1$.
Then
\begin{gather*}
L_n(x)=[a_{n}(x),b_{n}(x)]
\end{gather*}
and $L(x)$ is either the real line, a~left half-closed inf\/inite interval, or a~right half-closed inf\/inite interval.
To prove this claim notice that $\{a_{n}(x)\}_{n}$ is a~decreasing sequence and $\{b_{n}(x)\}_{n}$ is an increasing
sequence.
Indeed, if $x_{n}=0$ then $a_{n+1}(x)=a_{n}(x)$ and $b_{n}(x)<b_{n+1}(x)$, and if $x_{n}=1$ then $a_{n+1}(x)<a_{n}(x)$
and $b_{n}(x)=b_{n+1}(x)$.
It follows that if there is $n\in\mathbb{N}$ such that $x_{i}=0$ for all $i\ge n$, then $L(x)=[a_{n}(x),\infty)$ and, if
$x_{i}=1$ for all $i\ge n$, then $L(x)=(-\infty,b_{n}(x)]$.
Otherwise $L(x)=\mathbb{R}$.

We claim that for any $x\in X$ we have that $a_{n}(x)=-\sum\limits_{j=0}^{n-1}x_{j}2^{j}$ and
$b_{n}(x)=2^{n}-\sum\limits_{j=0}^{n-1}x_{j}2^{j}$ for all $n\ge1$.
In particular, $b_{n}(x)-a_{n}(x)=2^{n}$.
One can prove this claim by induction as follows.
Let $x\in X$.
If $x_{0}=0$ then $a_{1}(x)=F_{0}^{-1}(0)=0$ and $b_{1}(x)=F_{0}^{-1}(1)=2$.
If $x_{0}=1$ then $a_{1}(x)=F_{1}^{-1}(0)=-1$ and $b_{1}(x)=F_{1}^{-1}(1)=1$.
Thus $a_{1}(x)=-x_{0}2^{0}$ and $b_{1}(x)=2-x_{0}2^{0}$ and the claim holds for $n=1$.
Suppose that the induction hypothesis holds for all inf\/inite words for $n\ge1$.
Then, if $x\in X$ and $z:=x_{1}x_{2}\cdots$, it follows that
$a_{n}(z)=-\sum\limits_{j=0}^{n-1}z_{j}2^{j}=-\sum\limits_{j=1}^{n}x_{j}2^{j-1}$ and
$b_{n}(z)=2^{n}-\sum\limits_{j=0}^{n-1}z_{j}2^{j}=2^{n}-\sum\limits_{j=1}^{n}x_{j}2^{j-1}$.
If $x_{0}=0$ then
\begin{gather*}
a_{n+1}(x)=2a_{n}(z)=-\sum\limits_{j=1}^{n}x_{j}2^{j}=-\sum\limits_{j=0}^{n}x_{j}2^{j}
\end{gather*}
and
\begin{gather*}
b_{n+1}(x)=2b_{n}(z)=2^{n+1}-\sum\limits_{j=1}^{n}x_{j}2^{j}=2^{n+1}-\sum\limits_{j=0}^{n}x_{j}2^{j}.
\end{gather*}
If $x_{0}=1$ then
\begin{gather*}
a_{n+1}(x)=2a_{n}(z)-1=-1-\sum\limits_{j=1}^{n}x_{j}2^{j}=-\sum\limits_{j=0}^{n}x_{j}2^{j}
\end{gather*}
and
\begin{gather*}
b_{n+1}(x)=2b_{n}(z)-1=2^{n+1}-1-\sum\limits_{j=1}^{n}x_{j}2^{j}=2^{n+1}-\sum\limits_{j=1}^{n}x_{j}2^{j}.
\end{gather*}
Thus the induction holds and our claim is proved.

Therefore, for $n\in\mathbb{N}$ and $\alpha\in W^{n}$ we have that $L_{n}(\alpha)=[a_{n}(\alpha),b_{n}(\alpha)]$ and
$L_{n}=\bigsqcup\limits_{\vert\alpha\vert=n}Z(\alpha)\times [a_{n}(\alpha),b_{n}(\alpha)]$.
The local homeomorphism $\tilde{\sigma}$ on~$L$ is
\begin{gather*}
\tilde{\sigma}(x,t)=
\begin{cases}
(\sigma(x),t/2) & \text{if}\quad x_0=0,
\\
(\sigma(x),t/2+1/2) & \text{if}\quad x_0=1.
\end{cases}
\end{gather*}

Note that~$L$ is not locally compact.
For example, the point $(x,0)\in L$ does not have a~compact neighborhood, where~$x$ is the sequence for which $x_i=0$
for all $i\in\mathbb{N}$.
To see this, recall that if $C\subset L$ is compact then there is $n\ge 0$ such that $C\subset L_n$
(see~\cite[p.~2]{CurPat_Top88}).
However, one can check that there is no open set in~$L$ containing $(x,0)$ that is a~subset of any of the $L_n$'s.
\end{Example}

\begin{Example}
Fix $N > 1$, $r \in (0, 1)$ and let $e_1, \dots, e_N \in \mathbb{R}^N$ be the standard basis elements.
For $j = 1, \dots, N$, def\/ine $F_j: \mathbb{R}^N \to \mathbb{R}^N$ by $F_j(x) = rx + (1-r)e_j$.
Then $(F_1, \dots, F_N)$ forms an iterated function system on $Y = \mathbb{R}^N$ (endowed with the usual metric).
If $N=2$ and $r=1/2$, then the invariant set~$K$ is homeomorphic to the unit interval and $(F_1,F_2)$ is conjugate to
the iterated function system of the previous example.
If $N=3$ and $r=1/2$ then~$K$ is homeomorphic to the Sierpinski gasket.

We have $W = \{1, \dots, N \}$, $X = W^\infty = \{1, \dots, N \}^\infty$.
Note that for $x \in X$ we have
\begin{gather*}
\lim_{k\to\infty} F_{x_1} \circ \dots \circ F_{x_k}(0) = (1-r)\sum\limits_{j=1}^\infty r^{j-1}e_{x_j}.
\end{gather*}
Hence, the invariant set is given by $K = \Big\{(1-r)\sum\limits_{j=1}^\infty r^{j-1}e_{x_j}: x \in X \Big\}$.
It is straightforward to check that the iterated function system is totally disconnected if $r < 1/2$.

Observe that $F_j^{-1}(x) = \frac{1}{r}(x - (1-r)e_j)$.
For $\alpha \in W^n$, $L_{n}(\alpha)$ consists of all points of the form
\begin{gather*}
\frac{1-r}{r^{n+1}}\left(\sum\limits_{j=1}^\infty r^{j}e_{x_j} - \sum\limits_{i=1}^{n} r^{i}e_{\alpha_i}\right).
\end{gather*}
And, of course, for each $x \in X$, we have $L(x) = \bigcup_n L_n(x)$ (with the inductive limit topology).

We will show that the groupoid $\mathcal{G}\ast L$ is not minimal.

Let $y \in X$ be the sequence for which $y_i = 1$ for all~$i$.
We claim that the orbit of $(y, e_1) \in L$ is not dense.
In particular we show that $(y, e_2)$ is not in the closure of $\{\gamma\cdot(y, e_1): s(\gamma) = y \}$.
If $s(\gamma) = y$, then $\gamma = (x, m - n, y)$ (note $x_i = 1$ for $i \ge m$).
Since $e_1$ is a~f\/ixed point for $F_1$
\begin{gather*}
\gamma\cdot(y, e_1) = \big(x,F_{x(m)}^{-1}\circ F_{y(n)}(e_1)\big) = \big(x,F_{x(m)}^{-1}(e_1)\big)
\end{gather*}
Def\/ine $f: \mathbb{R}^N \to \mathbb{R}$ by $f(t_1, \dots, t_N) = t_2$.
Note that $F_j^{-1}(t) = \frac{1}{r}(t - (1-r)e_j)$, so if $f(t) \le 0$, then $f(F_j^{-1}(t)) \le 0$.
Observe that projection $\pi_2: L \to \mathbb{R}^N$ is continuous.
An easy induction argument now shows that
\begin{gather*}
(f \circ \pi_2)(\gamma\cdot(y, e_1)) = f\big(F_{x(m)}^{-1}(e_1)\big) \le 0.
\end{gather*}
Since $(f \circ \pi_2)(y, e_2) = 1$ (and $f \circ \pi_2$ is continuous), $(y, e_2)$ is not in the closure of the orbit
of $(y, e_1)$.

If $r<1/2$ the iterated function system is totally disconnected, and, hence, $F_i(K)$ is open in~$K$.
Thus, if $r<1/2$,~$L$ is locally compact.
However, if $r\ge 1/2$ then~$L$ is not locally compact.
An argument similar with the one at the end of Example~\ref{ex:unit} shows that, if $r\ge 1/2$, the point $(y,e_1)\in L$
does not have a~compact neighborhood in~$L$, where~$y$ is the sequence for which $y_i=1$ for all $i\in \mathbb{N}$.
\end{Example}

\subsection*{Acknowledgements}
The work of the f\/irst author was partially supported by a~grant from the Simons Foundation (\#209277 to Marius Ionescu).
The authors would like to thank the referees for their helpful comments.

\pdfbookmark[1]{References}{ref}
\LastPageEnding

\end{document}